\documentclass[11pt,a4paper]{article}
\usepackage{amsthm,amssymb,amsmath}
\usepackage{graphicx}
\newcommand{\R}{\mathbb R}
\newcommand{\henbi}[2]{\frac{\partial #1}{\partial #2}}
\newcommand{\pk}[1]{\mathrm{P}_{#1}}
\newcommand{\lagk}[1]{\Pi_h^{(#1)}}
\newcommand{\set}[1]{\left\{#1\right\}}
\newcommand{\rey}{\mathrm{Re}}
\newtheorem{theorem}{Theorem}{\bf}{}
\theoremstyle{definition}
\newtheorem*{Thschemezettai}{Scheme LG-LLV}{\bf}{}
\newtheorem{exam}{Example}{\bf}{}
\title{Numerical computation of triangular cavity flows by a Lagrange--Galerkin scheme with a locally linearized velocity}
\author{Masahisa Tabata\thanks{Department of Mathematics, Waseda University, Tokyo, Japan, tabata@waseda.jp}, ~
Shinya Uchiumi\thanks{Graduate School of Fundamental Science and Engineering, Waseda University, Tokyo, Japan, su48@fuji.waseda.jp} \thanks{Research Fellow of Japan Society for the Promotion of Science}
}
\date{}
\begin{document}
\maketitle

\begin{abstract}
We show numerical results of triangular cavity flow problems solved by a Lagrange--Galerkin scheme free from numerical quadrature.  
The scheme has recently developed by us, where a locally linearized velocity and the backward Euler approximation are used in finding the position of fluid particle at the previous time step.    
Since the scheme can be implemented exactly as it is, the theoretical stability and convergence results are assured, while the conventional Lagrange--Galerkin schemes may encounter the instability caused by numerical quadrature errors.  
The scheme is employed to solve cavity flow problems in triangular domains, where we observe the bifurcation of stationary solutions and the patterns of streamlines.

~

\end{abstract}

\section{Introduction}
The Lagrange--Galerkin method is a powerful numerical method for flow problems.  
It has such advantages that it is robust for convection-dominated problems and that the resultant matrix to be solved is symmetric.  
The Lagrange--Galerkin method is based on the characteristic method and it always contains the integration of composite function terms.   
Stability and convergence of the conventional Lagrange--Galerkin scheme have been proved under the condition that the integration is  computed exactly.
Since it is difficult to perform the exact integration in real problems, numerical quadrature is usually employed.    
It is, however, reported that instability may occur caused by numerical quadrature error in \cite{Tabata2007,TabataFujima2006,TabataUchiumi1,TabataUchiumiNS,TSTeng}.
Thus, it has been known to be a drawback of the Lagrange--Galerkin method that it may lose the stability when rough numerical quadrature is employed to integrate the composite function terms that characterize the method.     

Recently we have developed a Lagrange--Galerkin scheme with a locally linearized velocity \cite{TabataUchiumi1,TabataUchiumiNS}, which needs no numerical quadrature and overcomes the drawback of the instability mentioned above.
In this paper, after reviewing the scheme, we employ it to solve cavity flow problems in triangular domains and show the numerical results.  
  
The paper is organized as follows.  
In Section 2 we review the Lagrange--Galerkin scheme with a locally linearized velocity for the Navier--Stokes equations.  
In Section 3 we show numerical results solved by the scheme for cavity flow problems in triangular domains.  
The paper is concluded with some remarks.  
\section{A Lagrange--Galerkin scheme with a locally linearized velocity}
We review a Lagrange--Galerkin scheme with a locally linearized velocity \cite{TabataUchiumi1,TabataUchiumiNS}.  

Let us consider the Navier--Stokes problem:
find $(u,p):\Omega \times (0,T) \to \R^d \times \R$
such that
\begin{equation}\label{NS}
\begin{split}
	\henbi{u}{t}+(u \cdot \nabla) u
		- \nu \Delta u + \nabla p = f & \quad \text{in~} \Omega \times (0,T),\\
		\nabla \cdot u = 0 & \quad \text{in~} \Omega \times (0,T),\\
		u = 0 & \quad \text{on~} \partial \Omega \times (0,T),\\
		u = u^0 & \quad \text{in~} \Omega \text{~at~} t=0,
\end{split}
\end{equation}
where 
$\Omega$ is a polygonal or polyhedral domain of $\R^d ~ (d=2,3)$,
$\partial \Omega$ is the boundary of $\Omega$,
$T>0$ is a time and
$\nu>0$ is a viscosity. 
Functions $f:\Omega \times (0,T) \to \R^d$ 
and $u^0:\Omega \to \R^d$ are given.

Suppose the velocity field $u$ is smooth.
The characteristic curve $X(t; x,s)$ subject to the initial condition $X(s)=x$ 
is defined by the solution of the system of the ordinary differential equations, 
\begin{equation}
\begin{split}
\label{odes}
	\frac{dX}{dt}(t; x,s)&=u(X(t;x,s),t), \quad t<s,\\
	X(s;x,s)&=x.
\end{split}
\end{equation}
Then, we can write the material derivative term $(\henbi{}{t}+u\cdot \nabla)u$ at $(X(t),t)$ as follows: 
\begin{equation*}
\left(\henbi{u}{t} + (u \cdot \nabla) u \right)(X(t),t)=\frac{d}{dt} u(X(t),t).
\end{equation*}
Let $\Delta t>0$ be a time increment.  
Let 
$t^n \equiv n\Delta t$ and 
$\psi^n \equiv \psi(\cdot,t^n)$ for a function $\psi$ defined in $\Omega \times (0,T)$.
For $w:\Omega \to \R^d$ we define the mapping $X_1(w):\Omega \to \R^d$ by
\begin{equation*}\label{eq:x1def}
	(X_1(w))(x) \equiv x - w(x)\Delta t.
\end{equation*}
The image of $x$ by $X_1(u(\cdot,t))$ is nothing but the approximate value of $X(t-\Delta t;x,t)$ obtained by solving \eqref{odes} by the backward Euler method.
Then, it holds that
\begin{equation*}\label{eq:materialApprox}
\henbi{u^n}{t}+(u^n \cdot \nabla) u^n = \frac{u^n-u^{n-1}\circ X_1(u^{n-1})}{\Delta t} + O(\Delta t),
\end{equation*}
where $\circ$ stands for the composition of functions. 

Let $\mathcal T_h$ be a triangulation of $\bar \Omega$ and 
$h \equiv \max_{K\in \mathcal T_h} \operatorname{diam}(K)$ the maximum element size.
We always consider a regular family of triangulations $\{ \mathcal{T}_h  \}_{h\downarrow 0}$.  
Let $V_h \times Q_h$ be the $\pk2$/$\pk1$-finite element (or Hood--Taylor element) 
\begin{equation*}
\begin{split}
	V_h &\equiv \{v_h\in C(\bar{\Omega})^d; ~ v_{h|K}\in \pk 2(K)^d, \forall K \in \mathcal T_h, ~ v_{h|\partial \Omega}=0\}, \\
	Q_h &\equiv \left\{ q_h\in C(\bar{\Omega}); ~ q_{h|K}\in \pk 1(K), \forall K \in \mathcal T_h, ~ \int_{\Omega} q_h dx = 0 \right\},
\end{split}
\end{equation*}
where $C(\bar{\Omega})$ is the set of continuous functions on $\bar{\Omega}$ and
$\pk k (K)$ is the set of polynomials on $K$ whose degrees are less than or equal to $k$.
We denote by $\Pi_h^{(1)}$ the Lagrange interpolation operator to the $\pk 1$-finite element space.

Let $u_h^0$ be an approximation of $u^0$.
A Lagrange--Galerkin scheme with a locally linearized velocity \cite{TabataUchiumiNS} is as follows. 
\begin{Thschemezettai}
Find $\set{(u_h^n, p_h^n)}_{n=1}^{N_T} \subset V_h \times Q_h$ such that 
\begin{subequations}\label{eq:scheme0}
\begin{alignat*}{2}
	\left( \frac{u_h^{n} - u_h^{n-1}\circ  X_1(\lagk 1 u_h^{n-1})}{\Delta t}, v_h\right) + a(u_h^{n}, v_h)
	+ b(v_h, p_h^{n}) 	&= (f^{n},v_h), ~  	&&\forall v_h \in V_h, \\
	b(u_h^{n}, q_h) &= 0, ~ 	&&\forall q_h \in Q_h,
\end{alignat*}
\end{subequations}
for $n=1,\dots, N_T$,
where $N_T \equiv \lfloor T/\Delta t \rfloor$.
A pair of parentheses $(\cdot , \cdot)$ shows the $L^2(\Omega)$-inner product $(f,g) \equiv \int_{\Omega} f g ~ dx$. 
The inner products in $L^2(\Omega)^d$ and $L^2(\Omega)^{d\times d}$ are also denoted by the same notation.
The bilinear forms $a$ and $b$ are defined by 
\begin{equation*}
	a(u,v) \equiv \nu(\nabla u,\nabla v), \quad
	b(v,q) \equiv -(\nabla \cdot v,q).
\end{equation*}
\end{Thschemezettai}

\begin{figure}
	\centering
	\includegraphics[width=2in]{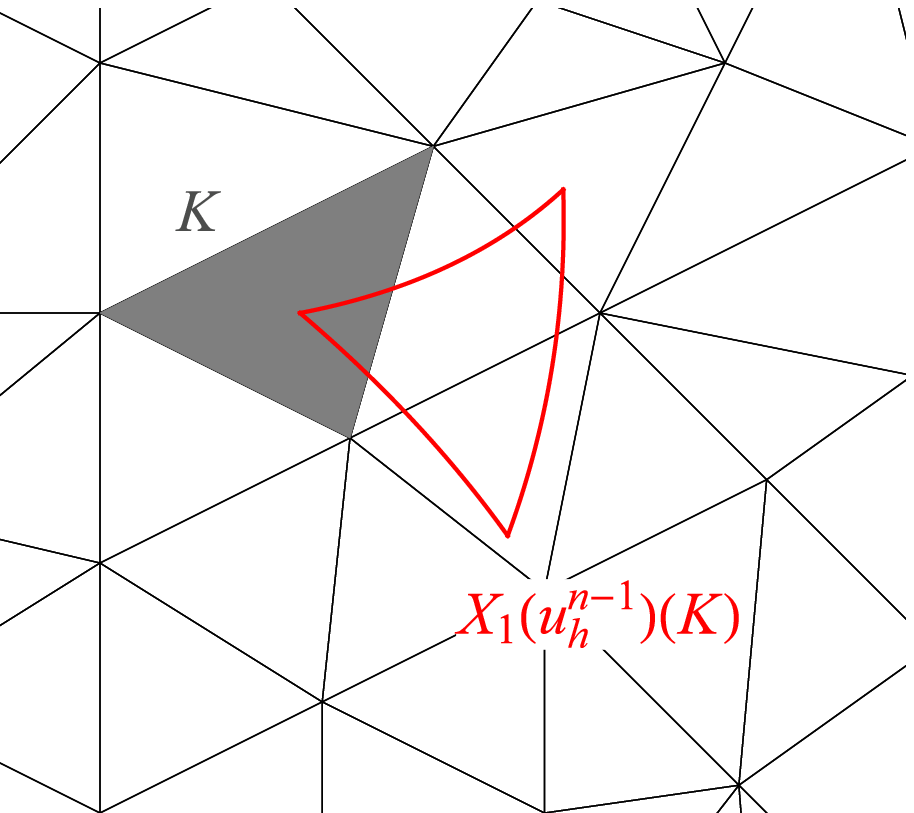} \quad
	\includegraphics[width=2in]{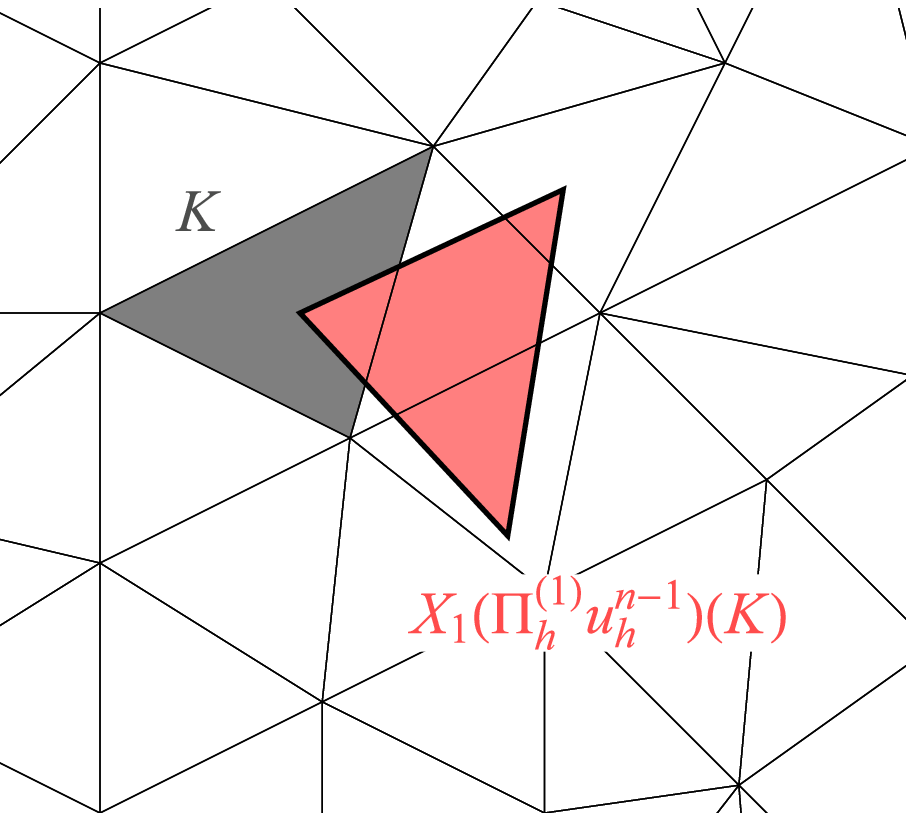}
	\caption{The image $X_1(u_h^{n-1})(K)$ (left) and $X_1(\Pi_h^{(1)} u_h^{n-1})(K)$ (right).}
	\label{fig:kx1u}
\end{figure}
In the conventional Lagrange--Galerkin scheme, the original velocity $u_h^{n-1}$ is used for the composite term.  
Since the function $u_h^{n-1}$ is quadratic, 
the image $X_1(u_h^{n-1})(K)$ is, in general, a curved triangle (Fig. \ref{fig:kx1u}, left).
Hence, it is hard to calculate the composite function term $( u_h^{n-1} \circ X_1(u_h^{n-1}), v_h )$ exactly.
In practice, the numerical quadrature has been used for the composite function term.
However, it has been reported that numerical quadrature causes the instability \cite{Tabata2007,TabataFujima2006,TabataUchiumi1,TabataUchiumiNS,TSTeng}.
In Scheme LG-LLV, thanks to the introduction of the locally linearized velocity, the image $X_1(\Pi_h^{(1)} u_h^{n-1})(K)$ becomes a triangle (Fig. \ref{fig:kx1u}, right), which makes the exact integration possible.

Let $u_h^0$ be the first component of the Stokes projection \cite{TabataUchiumiNS} of $(u^0,0)$.  
Suppose the family of triangulations $\{\mathcal T_h\}_{h \downarrow 0}$ satisfies the inverse assumption \cite{Ciarlet} and that 
$(u,p)$ is ``smooth'' (For the precise regularity see \cite{TabataUchiumiNS}).  
Then, we have the following convergence result \cite[Theorem 1]{TabataUchiumiNS}.  
\begin{theorem}\label{theo:mainConvergence}
Let 
$V_h\times Q_h$ be the $\pk2/\pk1$-finite element space.
Then, there exist positive constants  
$c_0$ and $h_0$ 
such that if $h\in (0,h_0]$ and $\Delta t \leq c_0 h^{d/4}$, 
the solution $(u_h,p_h) \equiv \set{(u_h^n,p_h^n)}_{n=0}^{N_T}$ of Scheme LG-LLV exists  
and the estimates
\begin{equation*}\label{eq:mainestimate}
	\| u_h-u \|_{\ell^\infty(H^1)}, \| p_h-p \|_{\ell^2(L^2)} \leq c_1 (h^2+\Delta t) 
\end{equation*}
hold, where $c_1$ is a positive constant independent of $h$ and $\Delta t$.
\end{theorem} 
Here we have used the norms defined by 
\begin{align*}
\| v_h \|_{\ell^{\infty}(H^1)}&=\max\{ \| v_h^n \|_{H^1(\Omega)}; ~n=0,\cdots,N_T \}, \\
\| v_h \|_{\ell^{2}(L^2)}&=\left\{ \Delta t \sum_{n=1}^{N_T}  \| v_h^n \|_{L^2(\Omega)}^2 \right\}^{1/2}.  
\end{align*}
By virtue of this theorem it is assured that the solution of LG-LLV, which is exactly computable, converges to the exact solution.  
\section{Numerical results}\label{sec:numex}
We solve cavity flow problems in triangular domains by Scheme LG-LLV.  
In the cavity flow problem a discontinuous boundary condition that the velocity $u=(1,0)$ on a side parallel to the $x_1$-axis and $u=0$ on the other boundary is often imposed, e.g., Erturk--Gokcol \cite{ErturkGokcol}.  
It is, however, known that under this condition there is no weak solution $(u,p)$ of \eqref{NS}.  
In order to assure the existence of the weak solution of \eqref{NS} we deal with regularized cavity flow problems, where the prescribed velocity is continuous on the boundary.
In the computation we use the following criterion to judge whether the stationary state is numerically attained or not, 
\begin{equation*}
	\max_P \frac{|\psi^n(P) - \psi^{n-1}(P)|}{\Delta t} < 10^{-4}, 
\end{equation*}
where $P$ runs all nodes and $\psi=u_h$ and $p_h$.  
\begin{figure}
\centering
	\includegraphics[width=2in]{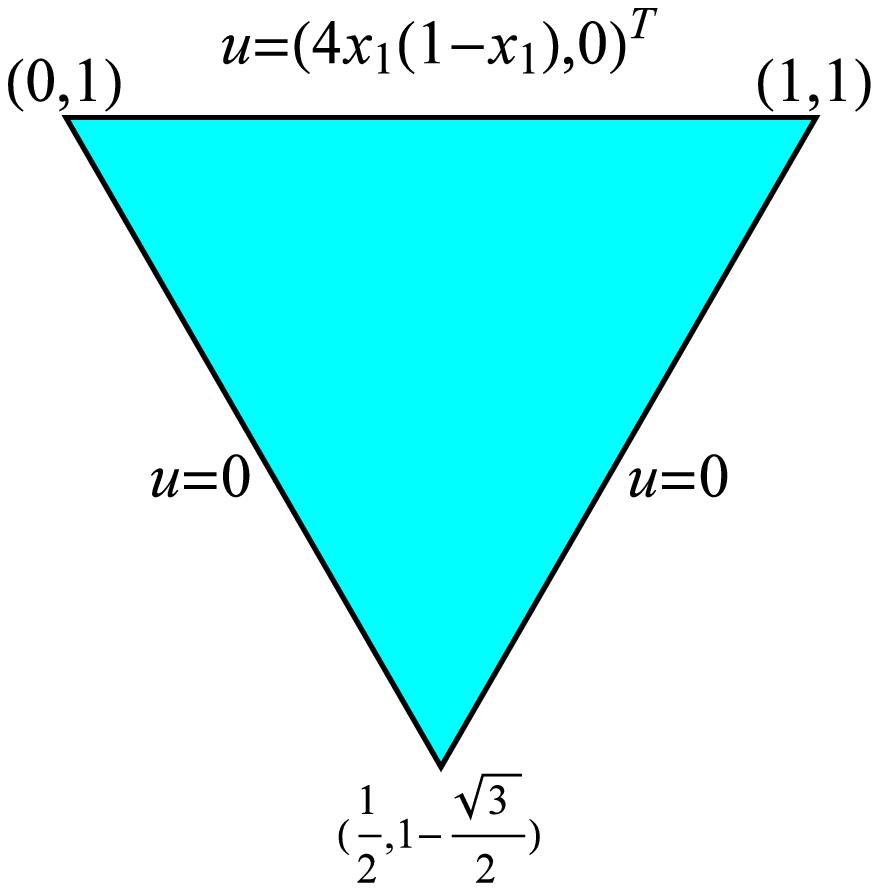}
	\includegraphics[width=2in]{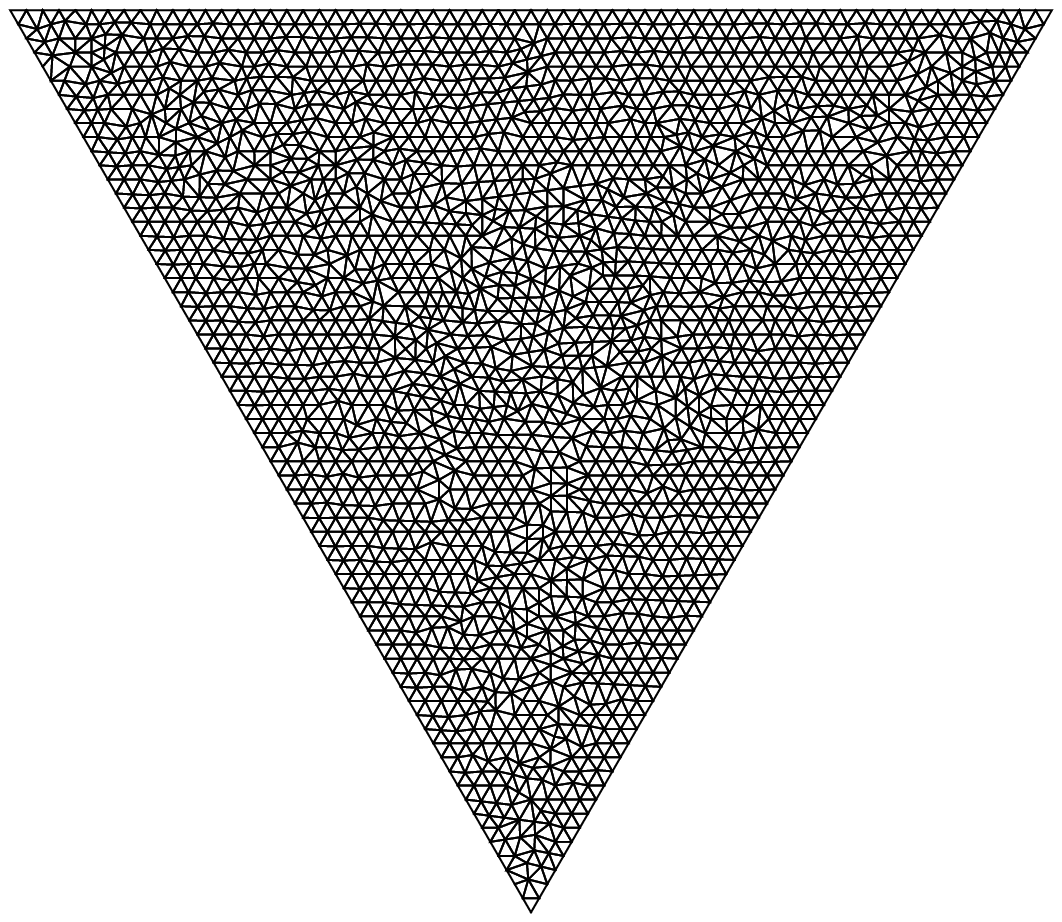}
	\caption{An equilateral triangular domain $\Omega$ and boundary conditions (left), and the triangulation of $\Omega$ (right).}
	\label{fig:domaintriAndtri64}
\end{figure}
\begin{exam}\label{exam:tri}
The domain $\Omega$ and boundary conditions are stated in Fig. \ref{fig:domaintriAndtri64} (left).
We set $\nu=1/\rey$ and $f=0$.
\end{exam}
\begin{figure}
\centering
	\includegraphics[width=2in]{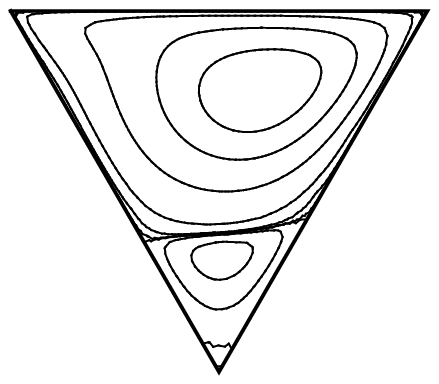}
	\includegraphics[width=2in]{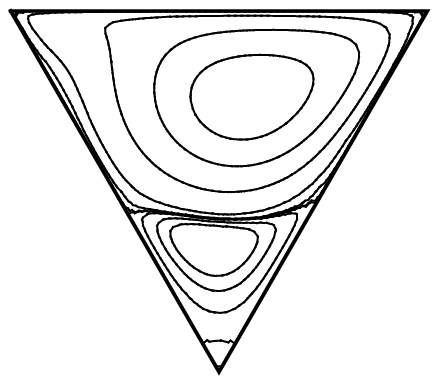}
	\\
	\includegraphics[width=2in]{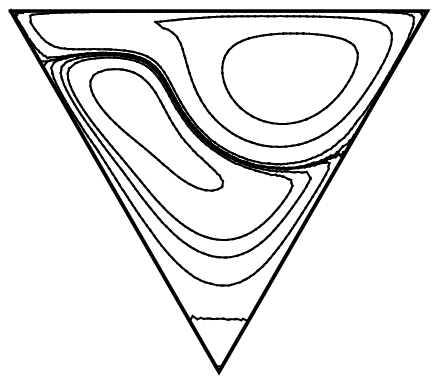}
	\includegraphics[width=2in]{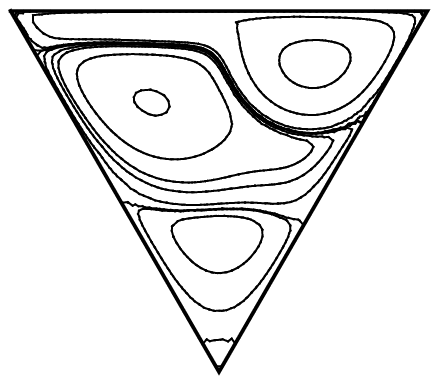}
	\caption{The streamlines of the stationary solutions $u_h^{*}$
	starting from $u_h^0=0$ 
	in Example \ref{exam:tri}. $\rey=500$ (top left), $1000$ (top right), $2000$ (bottom left) and $4000$ (bottom right).}
	\label{fig:cavtri}
\end{figure}
We use the mesh obtained by FreeFem++ \cite{FreeFemCite} by dividing each side into 64 segments, see Fig.  \ref{fig:domaintriAndtri64} (right).
We set the initial velocity $u_h^0$ to be zero and vary the Reynolds numbers, $\rey=500, 1000, 2000$ and $4000$.
The time increment is chosen as $\Delta t=1/64$.
Figure \ref{fig:cavtri} shows the streamlines of the numerically stationary solutions at these Reynolds numbers.
We denote by $u_h^{*}(M)$ the stationary solution at $\rey=M$.  
The secondary vortex of $u_h^{*}(2000)$ is much larger than that of $u_h^{*}(1000)$.
The change is larger than that in the well-known square cavity flow problem at the same Reynolds numbers.

\begin{figure}
\centering
	\includegraphics[width=2in]{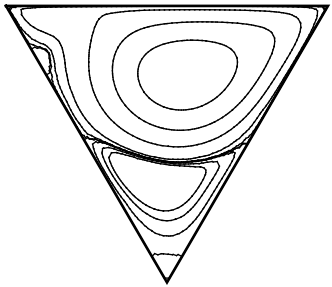}
	\caption{The streamlines of the stationary solution $u_h^{**}(2000)$
	in Example \ref{exam:tri}
	at $\rey=2000$
	starting from $u_h^0=u_h^{*}(1000)$.  
	}
	\label{fig:tricav_var}
\end{figure}
Next, we set the initial value $u_h^0$ to be $u_h^{*}(1000)$, the stationary solution at $\rey=1000$ obtained above.
The time increment is chosen as $\Delta t=1/256$.
Figure \ref{fig:tricav_var} shows the streamlines of the numerically stationary solution at $\rey=2000$. 
We denote it by $u_h^{**}(2000)$.  
It is observed that the solution is closer to the solution $u_h^{*}(1000)$ and is far different from $u_h^{*}(2000)$.
Here, we have used a smaller time increment $\Delta t=1/256$.
If we choose $\Delta t=1/64$, we 
get the solution $u_h^{*}(2000)$.
We also checked that, for $\Delta t=1/256$, we get again the same solution of $u_h^{*}(2000)$ starting from 
$u_{h}^0=0$.

\begin{figure}
	\centering
\includegraphics[width=1.9in]{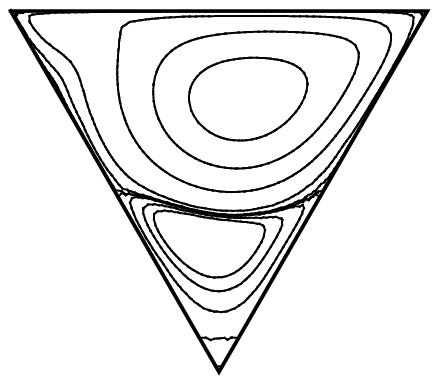}
\includegraphics[width=1.9in]{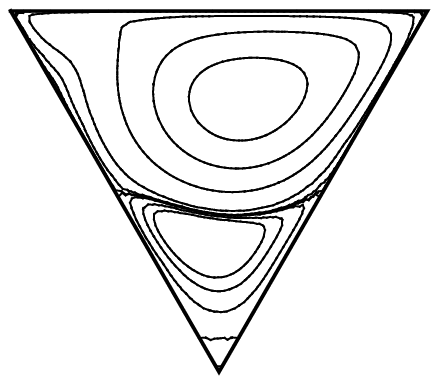} \\
\includegraphics[width=1.9in]{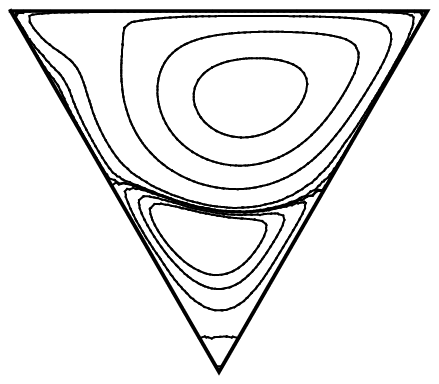}
\includegraphics[width=1.9in]{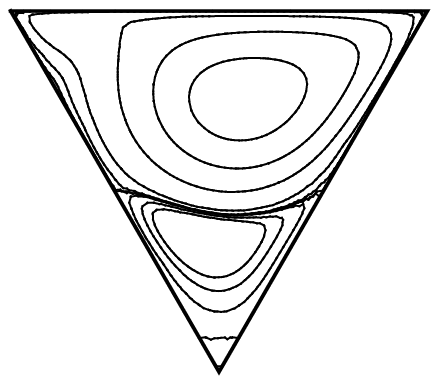} \\
\includegraphics[width=1.9in]{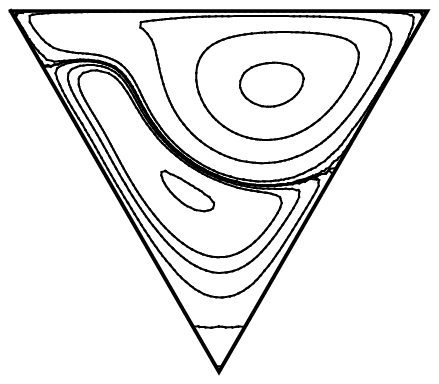}
\includegraphics[width=1.9in]{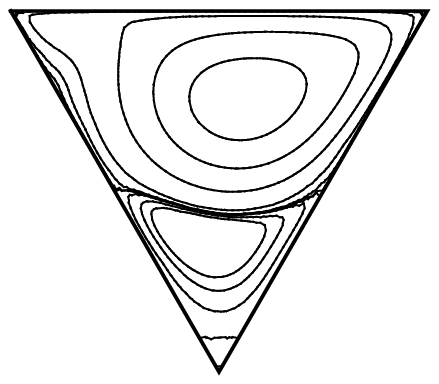} \\
\includegraphics[width=1.9in]{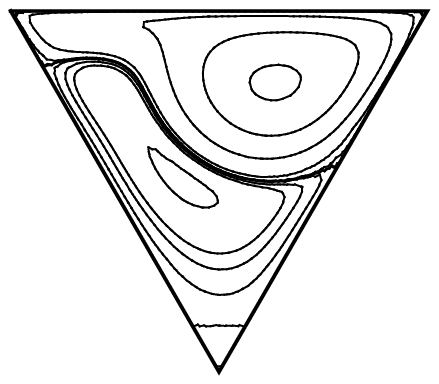}
\includegraphics[width=1.9in]{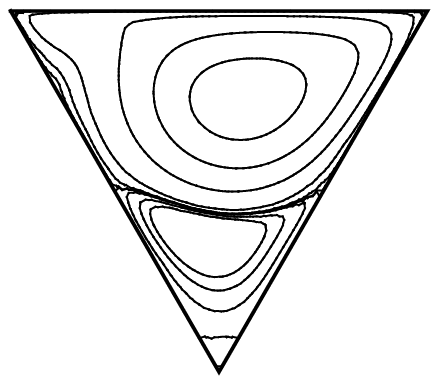} \\
\includegraphics[width=1.9in]{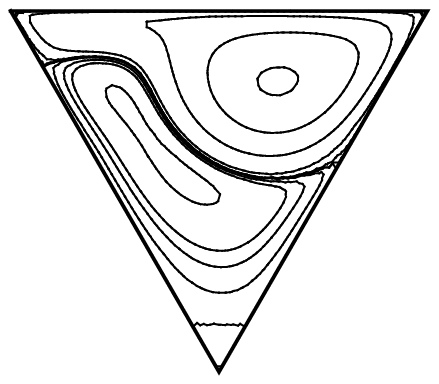}
\includegraphics[width=1.9in]{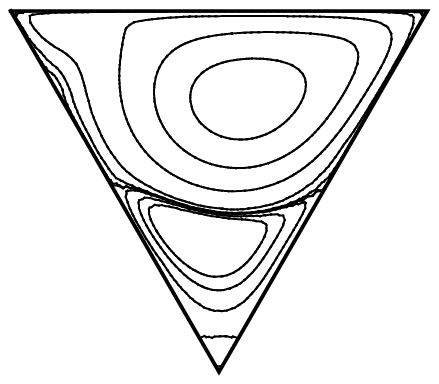} \\
	\caption{The streamlines of the stationary solutions of $u_h^{*}$ (left column) and of $u_h^{**}$ (right column) at $\rey=1550$, 1600, 1650, 1700, 1750 (from top to bottom).}
	\label{fig:10figs}
\end{figure}
\begin{figure}
	\centering
	\includegraphics[width=2in]{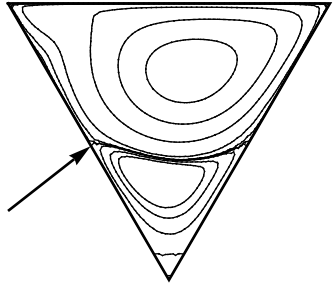}
	\caption{The stagnation point.}
	\label{fig:stag0}
\end{figure}
\begin{figure}
	\centering
	\includegraphics[width=4in]{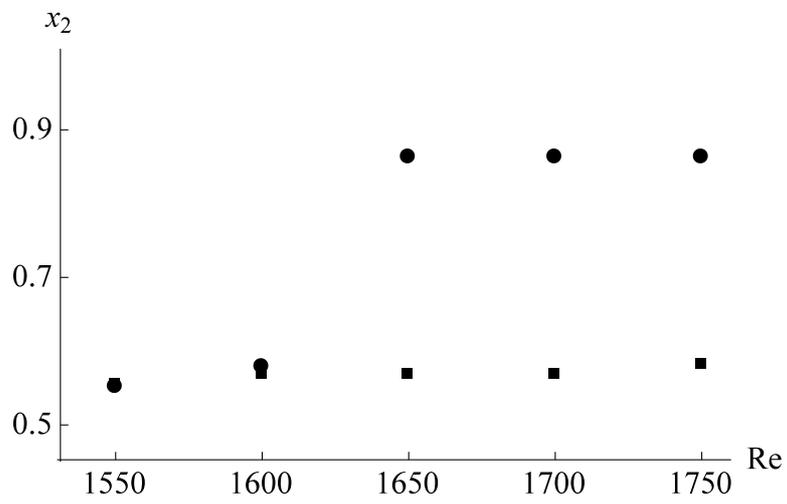}
	\caption{The $x_2$-coordinate of the stagnation point. 
		${\Large\bullet}$: the solutions $u_h^{*}$ and $\blacksquare$ : the solutions $u_h^{**}$.
	}
	\label{fig:stag}
\end{figure}
We vary the Reynolds numbers, $\rey=1500, 1550, 1600, 1650, 1700$ and $1750$.
We obtain stationary solutions at these Reynolds numbers in two ways. 
The first way is simple. 
We take the initial value $u_h^0=0$ and get the stationary solutions $u_h^{*}$ for all Reynolds numbers. 
The second way is as follows. 
We take the initial value $u_h^0=u_h^{*}(1500)$, solve the problem for $\rey=1550$ and get the stationary solution $u_h^{**}(1550)$, 
which is used as the initial value $u_h^0$ for $\rey=1600$.  
Repeating this procedure, we finally get $u_h^{**}(1750)$.
The time increment is chosen as $\Delta t=1/64$ in both two ways.  
In Fig. \ref{fig:10figs} we show the streamlines of the stationary solutions $u_h^{*}$ and $u_h^{**}$ obtained as above.  

We notice the stagnation point on the left side of the domain shown in Fig. \ref{fig:stag0}.   
Figure \ref{fig:stag} shows the $x_2$-coordinates of the stagnation points of the stationary solutions $u_h^{*}$ and $u_h^{**}$. 
In finding the stagnation point we have used a simple algorithm which enforces the point to be a node of the mesh.   
Hence, this is a discrete graph of mesh width $h_0=1/64$.  
The $x_2$-coordinates of two solutions are same at $\rey=1550$. 
The $x_2$-coordinates take slightly different values at $\rey=1600$.
We observe that the $x_2$-coordinate of the solution $u_h^{*}$ increases sharply from $\rey=1600$ to $1650$ while 
that of the solution $u_h^{**}$ remains almost same.  
The $x_2$-coordinates of $u_h^{*}$ at $\rey=1650, 1700$, and $1750$ are same and those of $u_h^{**}$ at $\rey=1600, 1650$, and $1700$ are same, which means that the differences are less than $h_0$.  
From the observation above it is concluded that around $\rey=1600$ the bifurcation of the solutions occurs.  
In order to obtain the critical Reynolds number more accurate computation will be required.  

Erturk and Gokcol \cite{ErturkGokcol} have performed numerical computation in fine grids of cavity flow problems in triangular domains.  
Although their boundary condition is discontinuous, 
our solutions $u_h^{**}$ in Figs. \ref{fig:tricav_var} and \ref{fig:10figs} have good agreements with the 
corresponding solutions in \cite[Fig. 2] {ErturkGokcol}.  
Note that our Reynolds numbers must be divided by $2\sqrt 3$ to correspond to their Reynold numbers 
because our domain is smaller than theirs by $2\sqrt{3}$.  
In \cite{ErturkGokcol} there is no figure corresponding to the solutions $u_h^{*}$.  

\begin{figure}
	\centering
	\includegraphics[width=2in]{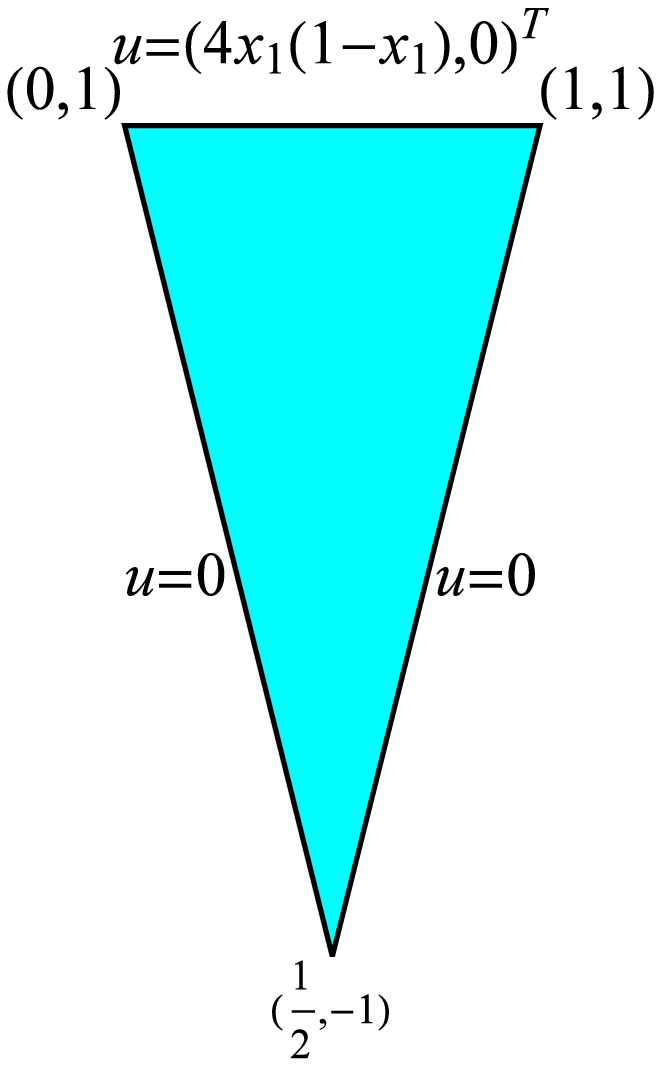}
	\includegraphics[width=2in]{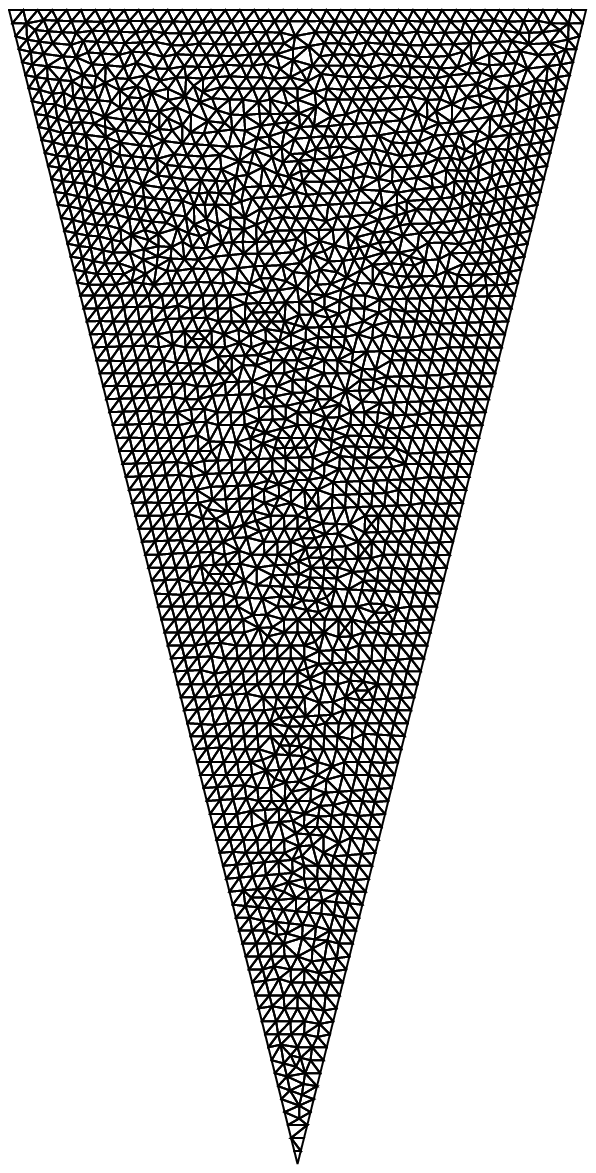}
	\caption{An isosceles triangular domain $\Omega$ and boundary conditions (left), and the triangulation of $\Omega$ (right).}
	\label{fig:domaintriAndtri2}
\end{figure}

\begin{exam}\label{exam:tri2}
The domain $\Omega$ and boundary conditions are stated in Fig. \ref{fig:domaintriAndtri2} (left).
We set $\nu=1/\rey$, $\rey=200, 400$, and $f=0$.
\end{exam}
\begin{figure}
	\centering
	\includegraphics[width=2in]{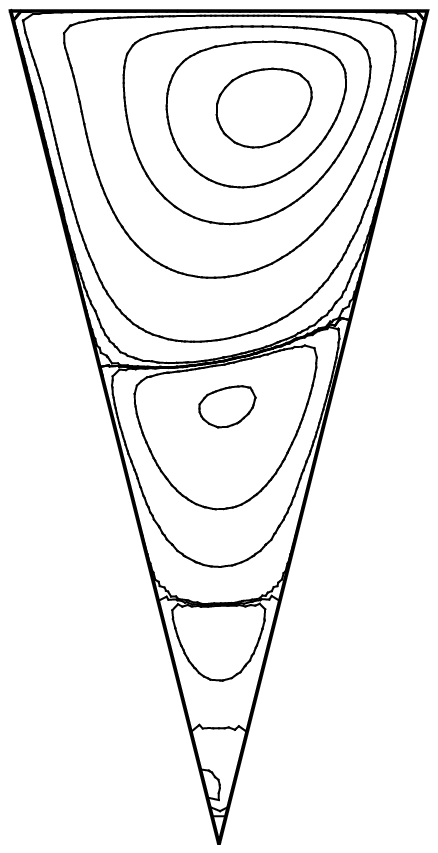}
	\includegraphics[width=2in]{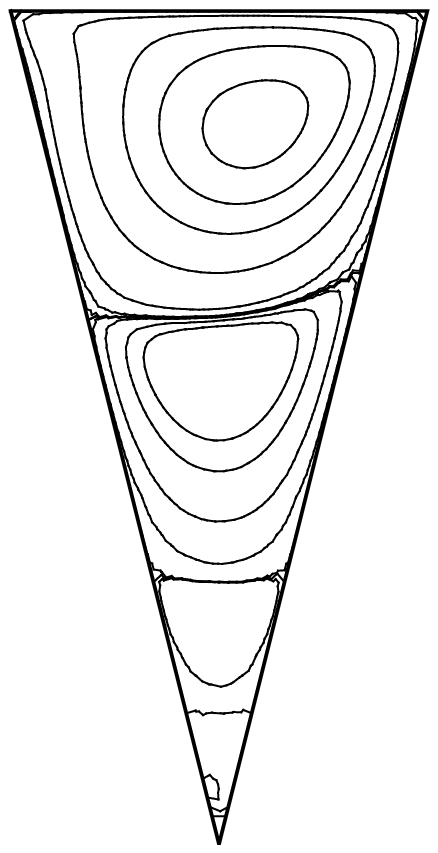}
	\caption{The streamlines of stationary solutions in Example \ref{exam:tri2}. 
		$\rey=200$  (left) and $400$ (right).
	}
	\label{fig:cavtri2}
\end{figure}
We use the mesh in Fig. \ref{fig:domaintriAndtri2} (right).
We show the streamlines of the stationary solutions starting from $u_h=0$ in Fig. \ref{fig:cavtri2}.
Although Erturk and Gokcol used the discontinuous boundary condition, theses figures look similar to their corresponding results in \cite[Fig. 6]{ErturkGokcol}.   
Note again that our Reynolds numbers must be divided by $2$ to correspond to their Reynold numbers 
because our domain is smaller than theirs by $2$.  

\section{Concluding remarks}\label{sec:conclusion}
We have shown numerical results of cavity flow problems in triangular domains by a Lagrange--Galerkin scheme free from numerical quadrature.
By virtue of the introduction of a locally linearized velocity, the scheme can be implemented exactly and the theoretical stability results are assured for practical numerical solutions.  
In the equilateral triangular domain, we have observed the bifurcation of the stationary solutions.  
We also got the stationary solutions at $\rey=200$ and 400 in an isosceles triangular domain, which show similar patterns of streamlines to the previous results \cite{ErturkGokcol}.  
We have solved the cavity flow problems subject to a continuous boundary condition to ensure the existence of the weak solution.  
To find numerically the critical Reynolds number of the bifurcation for such regularized cavity flow problems will be an interesting problem.   

\section*{Acknowledgment}
This work was supported by JSPS (the Japan Society for the Promotion of Science) under the Japanese-German Graduate Externship (Mathematical Fluid Dynamics).
The first author was supported by JSPS 
under Grants-in-Aid for Scientific Research (C), No. 25400212 and (S), No. 24224004 
and by Waseda University under Project research, Spectral analysis and its application to the stability theory of the Navier-Stokes equations of Research Institute for Science and Engineering.
The second author was supported by JSPS under Grant-in-Aid for JSPS Fellows, No. 26$\cdot$964.

\end{document}